\documentclass[twoside,leqno,twocolumn]{article}
\usepackage{ltexpprt}

\usepackage{bm}
\usepackage{amsfonts,amssymb,amsmath,dsfont} 
\usepackage{enumitem}
\setlist{noitemsep,topsep=0pt,parsep=0pt,partopsep=0pt}
\usepackage{xspace}
\usepackage{graphicx}
\usepackage{xcolor}
\usepackage{cite}
\usepackage{hyperref}
\usepackage{subfig}
\usepackage{comment}
\definecolor{darkgreen}{rgb}{0,0.4,0}
\definecolor{BrickRed}{rgb}{0.65,0.08,0}
\hypersetup{colorlinks=true,linkcolor=blue,citecolor=red,filecolor=BrickRed,urlcolor=darkgreen}
\usepackage{tikz}
\usepackage{tkz-tab}

\newcommand{\proba}{\mathds{P}}
\newcommand{\bigO}{\mathcal{O}}
\newcommand{\Pc}{\mathcal{P}}

\newcommand{\integers}{\mathds{Z}}

\newcommand{\naturals}{\integers_{\geq 0}}

\newcommand*{\eg}{\textit{e.g.,}\@\xspace}
\newcommand*{\ie}{\textit{i.e.,}\@\xspace}
			
\newcommand*{\resp}{\textit{resp.}\@\xspace}

\newcommand{\one}{\mathds{1}}
\newcommand{\id}{\operatorname{Id}}

\newcommand{\reach}{\operatorname{reach}}

\newcommand{\bridge}{\operatorname{Bridge}}

\newcommand{\OEIS}[1]{\href{http://oeis.org/#1}{OEIS~#1}}

\newcommand{\wo}{y} 
\newcommand{\wt}[1]{p_{#1}}
\newcommand{\Bb}{B_2}

\newcommand{\Z}{\mathbb{Z}}

\makeatletter
\newcommand{\I}{(1)}
\newcommand{\II}{(2)}
\newcommand{\III}{(3)}
\makeatother

\newtheorem{example}{Example}
\numberwithin{example}{section}

\newcommand{\proofofref}{}
\newproof{zproofof}{\proofofref}
\newenvironment{proofof}[1]
 {\renewcommand{\proofofref}{#1}\zproofof}
 {\endzproofof}

\begin{document}

\title{\Large Combinatorics of nondeterministic walks of the Dyck and Motzkin type\thanks{This paper will be published in the proceedings of ANALCO~2019.}}
\author{\'Elie de Panafieu\thanks{Nokia Bell Labs and Lincs, France} \\
\and
Mohamed Lamine Lamali\thanks{LaBRI \-- Universit\'e  de Bordeaux, France} \\
\and
Michael Wallner\footnotemark[3]}
\date{}

\maketitle







\begin{abstract} \small\baselineskip=9pt This paper introduces \emph{nondeterministic walks}, a new variant of one-dimensional discrete walks. At each step, a nondeterministic walk draws a random set of steps from a predefined set of sets and explores all possible extensions in parallel. 
We introduce our new model on Dyck steps with the nondeterministic step set $\{ \{-1\}, \{1\}, \{-1,1\}\}$ and Motzkin steps with the nondeterministic step set $\{\{-1\}, \{0\}, \{1\}, \{-1,0\}, \{-1,1\}, \{0,1\}, \{-1,0,1\}\}$.
For general lists of step sets and a given length, we express the generating function of nondeterministic walks where at least one of the walks explored in parallel is a bridge (ends at the origin). In the particular cases of Dyck and Motzkin steps, we also compute the asymptotic probability that at least one of those parallel walks is a meander (stays nonnegative) or an excursion (stays nonnegative and ends at the origin). 

This research is motivated by the study of networks involving encapsulations and decapsulations of protocols. Our results are obtained using generating functions and analytic combinatorics.

\textbf{Keywords.} Random walks, analytic combinatorics, generating functions, networking, encapsulation.
\end{abstract}

		\section{Introduction}

In recent years lattice paths have received a lot of attention in different fields, such as probability theory, computer science, biology, chemistry, physics, and much more~\cite{Kn69,HOP17,BP08}.
One reason for that is their versatility as models like \eg the up-to-date model of certain polymers in chemistry \cite{vRPR08}.
In this paper we introduce yet another application: 
the encapsulation of protocols over networks. 
To achieve this goal we generalize the class of lattice paths to so called \emph{nondeterministic lattice paths}.

	   \subsection{Definitions}

  \paragraph{Classical walks.}

\begin{figure}[ht]
	\centering
	\subfloat[\scriptsize A classical walk.]{%
    \resizebox{0.32\textwidth}{!}{\begin{tikzpicture}[shorten >=-3pt,shorten <=-3pt, x=1cm, y=0.5cm]
	\draw [very thin, gray,ystep=1] (-1,-2) grid (5,6) ;

	\draw [-*, very thick,red] (0,0) edge (1,2) ;
	
	\draw [-*, very thick, red] (1,2) edge (2,1) ;
	
	\draw [-*, very thick, red] (3,1) edge  (4,2) ;
	\draw [-*, very thick, red] (2,1)  edge (3,1) ;
	
	\draw [o->,thick] (-1,0) edge (4.9,0);
	\draw [o->,thick] (0,-1) edge (0,5.9);

	\foreach \x in {-2,...,6}
	\draw (-1,\x)--(-1,\x)  node[left] {\x};
	
	\foreach \y in {0,...,5}
	\draw (\y,-2)--(\y,-2)  node[below] {\y};
	\draw [-*, very thick, red] (0,0) ;	
	\end{tikzpicture}}
		\label{fig:walk}}
        \hfill
\subfloat[\scriptsize An N-walk.]{%
    \resizebox{0.32\textwidth}{!}{\begin{tikzpicture}[shorten >=-3pt,shorten <=-3pt, x=1cm, y=0.5cm]
	\draw [very thin, gray,ystep=1] (-1,-2) grid (5,6) ;

	\draw [-*, very thick] (0,0) edge[red] (1,2) (1,2) edge (2,3) (2,3) edge (3,3) edge[-] (3,1) (3,3) edge (4,5) edge (4,4) edge (4,3);
	
	\draw [-*, very thick, red] (1,2) edge (2,1) ;
	
	\draw [-*, very thick] (3,1) edge[-] (4,1) edge[red] (4,2) edge[-] (4,3) ;
	\draw [-*, very thick] (2,1)  edge (3,-1) edge[red] (3,1) (3,-1) edge (4,-1) edge (4,0) edge (4,1);
	
	\draw [o->,thick] (-1,0) edge (4.9,0);
	\draw [o->,thick] (0,-1) edge (0,5.9);

\foreach \x in {-2,...,6}
	\draw (-1,\x)--(-1,\x)  node[left] {\x};
	
\foreach \y in {0,...,5}
\draw (\y,-2)--(\y,-2)  node[below] {\y};
	\draw [-*, very thick, red] (0,0) ;	
	\draw [-*, very thick, red] (2,1) ;	
	\draw [-*, very thick, red] (3,1) ;	
	\end{tikzpicture}}
		\label{fig:n_walk}}
        \hfill
	\subfloat[\scriptsize Meanders compatible with the N-walk.]{%
    \resizebox{0.32\textwidth}{!}{\begin{tikzpicture}[shorten >=-3pt,shorten <=-3pt, x=1cm, y=0.5cm]
	\draw [very thin, gray,ystep=1] (-1,-2) grid (5,6) ;
	
	\draw [-*, very thick] (0,0) edge[red] (1,2) (1,2) edge (2,3) (2,3) edge (3,3) edge[-] (3,1) (3,3) edge (4,5) edge (4,4) edge (4,3);
	
	\draw [-*, very thick, red] (1,2) edge (2,1) (2,1) edge (3,1) ;
	
	\draw [-*, very thick] (3,1) edge (4,1) edge[red] (4,2) edge (4,3) ;

	\draw [o->,thick] (-1,0) edge (4.9,0);
	\draw [o->,thick] (0,-1) edge (0,5.9);

	\foreach \x in {-2,...,6}
	\draw (-1,\x)--(-1,\x)  node[left] {\x};
	
	\foreach \y in {0,...,5}
	\draw (\y,-2)--(\y,-2)  node[below] {\y};
	\draw [-*, very thick, red] (0,0) ;	
	\draw [-*, very thick, red] (2,1) ;	
	\draw [-*, very thick, red] (3,1) ;	
	\end{tikzpicture}}
		\label{fig:n_meander}} 
    \caption{Geometric realization of a walk, an N-walk, and its compatible meanders.}
   \label{fig:walks}
\end{figure}
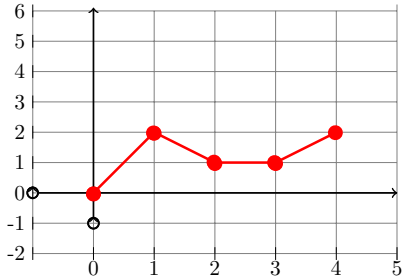
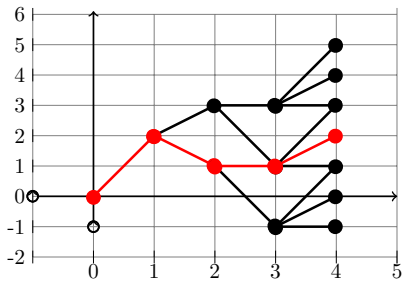
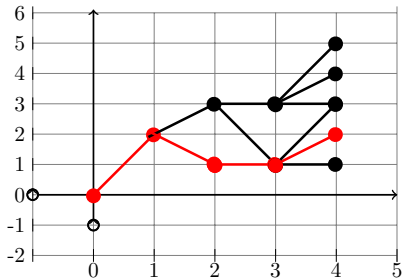

We mostly follow terminology from~\cite{BaFl02}.
%
Given a set $S$ 
of integers, called the \emph{steps}, 
a \emph{walk} is a sequence
$v = (v_1, \ldots, v_n)$ of steps $v_i \in S$. 
In this paper we will always assume that our walks start at the origin.
Its \emph{length} $|v|$ is the number $n$ of its steps,
and its \emph{endpoint} is equal to
the sum of the steps $\sum_{i=1}^n v_i$.
As illustrated in Figure~\ref{fig:walk},
a walk can be visualized by its \emph{geometric realization}.
Starting from the origin, the steps are added one by one to the previous endpoints. 
This gives a sequence $(\wo_j)_{0 \leq j \leq n}$ of ordinates at discrete time steps, such that
$\wo_0 = 0$ and $\wo_j := \sum_{i=1}^j v_i.$
A \emph{bridge} is a walk with endpoint $\wo_n=0$.
A \emph{meander} is a walk 
where all points have nonnegative ordinate,
\ie $\wo_j \geq 0$ for all $j=0,\ldots,n$.
An \emph{excursion} is a meander with endpoint $\wo_n=0$.

  \paragraph{Nondeterministic walks.}

This paper investigates a new variant of walks,
called \emph{nondeterministic walks}, or \emph{N-walks}.
In our context, this word does not mean ``\emph{random}''.
Instead it is understood in the same sense
as for automata and Turing machines.
A process is nondeterministic
if several branches are explored in parallel,
and the process is said to end in an accepting state
if one of those branches ends in an accepting state.
Let us now give a precise definition of these walks.
%


\begin{Definition}[Nondeterministic walks]
An \emph{N-step} is a nonempty set of integers.
Given an N-step set~$S$, 
an \emph{N-walk} $w$ 
is a sequence of N-steps.
Its \emph{length} $|w|$ is equal to the number of its N-steps.
\end{Definition}

\noindent As for classical walks we always assume that they start at the origin and we distinguish different types.

\begin{Definition}[Types of N-walks] \label{def:types_of_n_walks}
An N-walk $w = (w_1, \ldots, w_n)$ and a classical walk $v = (v_1, \ldots, v_n)$
are \emph{compatible} if they have the same length $n$, the same starting point, 
and for each $1 \leq i \leq n$,
the $i^{\text{th}}$ step is included in the $i^{\text{th}}$ N-step, \ie $v_i \in w_i$.
An \emph{N-bridge} (\resp \emph{N-meander}, \resp \emph{N-excursion})
is an N-walk compatible with at least one bridge
(\resp meander, \resp excursion).
Thus, N-excursions are particular cases of N-meanders.
\end{Definition}

\noindent The endpoints of classical walks are central to the analysis.
We define their nondeterministic analogues. 

\begin{Definition}[Reachable points] \label{def:reachable_points}
The \emph{reachable points} of a general N-walk are the endpoints
of all walks compatible with it.
For N-meanders, the reachable points are defined
as the set of endpoints of compatible meanders.
In particular, all reachable endpoints of an N-meander are nonnegative.
The minimum (\resp maximum) reachable point of an N-walk $w$ is denoted by $\min(w)$ (\resp $\max(w)$).
The minimum (\resp maximum) reachable point of an N-meander $w$ is denoted by $\min^+(w)$ (\resp $\max^+(w)$).
\end{Definition}
The geometric realization of an N-walk
is the sequence, for $j$ from $0$ to $n$,
of its reachable points after $j$ steps.
Figure~\ref{fig:walks} illustrates the geometric realization
of a walk $v = (2,-1,0,1)$ in (\ref{fig:walk}),
of an N-walk $w = \left(\{2\}, \{-1,1\}, \{-2,0\}, \{0,1,2\}\right)$ in (\ref{fig:n_walk}),
and of the classical meanders compatible with $w$ in (\ref{fig:n_meander}).
Note that the walk $v$ (highlighted in red) is compatible with the N-walk $w$.

  \paragraph{Probabilities.}

Any set of weights, and in particular any probability distribution
on the set of steps or N-steps
induces a probability distribution on walks or N-walks.
The probability associated to the walk or N-walk $w = (w_1, \ldots, w_n)$
is then the product $\prod_{i=1}^n \proba(w_i)$
of the probabilities of its steps or N-steps.

\subsection{Main results}


Our main results are the analysis of the asymptotic number of nondeterministic walks of the Dyck and Motzkin type with step sets $\left\{ \{-1\}, \{1\}, \{-1, 1\} \right\}$ and $\left\{\{-1\}, \{0\}, \{1\}, \{-1,0\}, \{-1,1\}, \{0,1\}, \{-1,0,1\} \right\}$, respectively. The results for the unweighted case where all weights are set equal to one are summarized in Table~\ref{tab:compareDyckMotzkin}.
These results are derived using generating functions and singularity analysis.
The reappearing phenomenon is the one of a simple dominating polar singularity
arising from the large redundancy in the steps.
The type of N-walk only influences the constant or the proportion among all N-walks. 
The lower order terms are exponentially smaller and of the square root type.
These are much more influenced by the types.
From a combinatorial point of view, we see a quite different behavior compared with classical paths. In particular, the limit probabilities for a Dyck N-walk of even length to be an N-bridge, an N-meander, or an N-excursion, are $1$, $1/2$, or $1/4$, and for Motzkin N-walks $1$, $3/4$, or $9/16$.


We also explore general N-steps and prove that the generating function of N-bridges is always algebraic.
N-excursions with general N-steps will be investigated in a longer version of this article.

\begin{table*}[ht]
	\begin{center}
	\begin{tabular}{|c||c|c|}
		\hline    Type  & Dyck N-steps  & Motzkin N-steps\\
        & 
        $\Pc(\{-1,1\}) \setminus \emptyset $ &
        $\Pc(\{-1,0,1\}) \setminus \emptyset $\\        
		\hline  
			N-Walk & 
            $3^n$ & 
            $7^n$ 
            \\
			N-Bridge & 
            $\frac{1+(-1)^n}{2} \left( 3^n - \frac{2\sqrt{2}}{\sqrt{\pi}} \frac{8^{n/2}}{\sqrt{n}} + \bigO\left(\frac{8^{n/2}}{n^{3/2}}\right) \right)$ & 
            $7^n - \sqrt{\frac{3}{\pi}} \frac{6^n}{\sqrt{n}} + \bigO\left(\frac{6^n}{n^{3/2}}\right)$ 
            \\
            N-Meander & 
            $\frac{3^n}{2} + \frac{3\sqrt{2}(1+(-1)^n) + 4(1-(-1)^n)}{\sqrt{\pi}} \frac{8^{n/2}}{\sqrt{ n^3}} + \bigO\left(\frac{8^{n/2}}{n^{5/2}}\right)$ & 
            $\frac{3}{4}7^n + \frac{3\sqrt{3}}{2 \sqrt{\pi}} \frac{6^n}{\sqrt{n^3}} + \bigO\left(\frac{6^n}{n^{5/2}}\right)$
            \\
            N-Excursion & $\frac{1+(-1)^n}{2}\left( \frac{3^n}{4} + 4\sqrt{2} \frac{8^{n/2}}{\sqrt{\pi n^3}}  + \bigO\left(\frac{8^{n/2}}{n^{5/2}}\right) \right)$ & 
            $\frac{9}{16}7^n - \gamma \frac{6^n}{\sqrt{\pi n^3}} + \bigO\left(\frac{6^n}{n^{5/2}}\right)$
            \\
		\hline
  \end{tabular}
\end{center}
\caption{The asymptotic number of nondeterministic unweighted (all weights equal to $1$) Dyck and Motzkin N-walks
with $n$ steps obeying different constraints: 
N-Bridges contain at least one classical bridge ending at $0$, N-Meanders contain at least one classical meander staying nonnegative, and N-excursions contain at least one classical excursion staying nonnegative and ending at $0$.
The constant $\gamma \approx 0.6183$ is an algebraic number defined as the positive real solution of $1024\gamma^4-8019\gamma^2+2916=0$.}
\label{tab:compareDyckMotzkin}
\end{table*}

	\subsection{Motivation and related work}

Let us start with a vivid motivation of the model using Russian dolls.
%
Suppose we have a set of $n+1$ people arranged in a line. There are three kinds of people. A person of the first kind is only able to put a received doll in a bigger one. A person of the second kind is only able to extract a smaller doll (if any) from a bigger one. If she receives the smallest doll, then she throws it away. Finally, a person of the third kind can either put a doll in a bigger one or extract a smaller doll if any. 
We want to know if it is possible for the last person to receive the smallest doll after it has been given to the first person and then, consecutively, handed from person to person while performing their respective operations. 
This is equivalent to asking if a given N-walk with each N-step $\in \left\{\{1\},\{-1\},\{-1,1\}\right\}$ is an N-excursion, \ie if the N-walk is compatible with at least one excursion. The probabilistic version of this question is: what is the probability that the last person can receive the smallest doll according to some distribution on the set of people over the three kinds?

	\paragraph{Networks and encapsulations.}

The original motivation of this work comes from networking. In a network, some nodes are able to encapsulate protocols (put a packet of a protocol inside a packet of another one), decapsulate protocols (extract a nested packet from another one), or perform any of these two operations (albeit most nodes are only able to transmit packets as they receive them). Typically, a tunnel is a subpath starting with an encapsulation and ending with the corresponding decapsulation. Tunnels are very useful for achieving several goals in networking (\eg interoperability: connecting IPv6 networks across IPv4 ones~\cite{wu2013transition}; security and privacy: securing IP connections~\cite{seo2005security}, establishing Virtual Private Networks~\cite{rosen2015multicast}, etc.). Moreover, tunnels can be nested to achieve several goals. Replacing the Russian dolls by packets, it is easy to see that an encapsulation can be modeled by a $\{1\}$ step and a decapsulation by a $\{-1\}$, while a passive transmission of a packet is modeled by a $\{0\}$ step.

Given a network with some nodes that are able to encapsulate or decapsulate protocols, a path from a sender to a receiver is \textit{feasible} if it allows the latter to retrieve a packet exactly as dispatched by the sender. 
Computing the shortest feasible path between two nodes is polynomial~\cite{LFCP16} if cycles are allowed without restriction. In contrast, the problem is  $\mathsf{NP}$-hard if cycles are forbidden or arbitrarily limited. In~\cite{LFCP16}, the algorithms are compared through worst-case complexity analysis and simulation. The simulation methodology for a fixed network topology is to make encapsulation (\resp decapsulation) capabilities available with some probability $p$ and observe the processing time of the different algorithms. It would be interesting, for simulation purposes, to generate random networks with a given probability of existence of a feasible path between two nodes. This work is the first step towards achieving this goal, since our results give the probability that any path is feasible (\ie is a N-excursion) according to a probability distribution of encapsulation and decapsulation capabilities over the nodes.

	\paragraph{Lattice paths.}

Nondeterministic walks naturally connect between lattice paths and branching processes. 
This is underlined by our usage of many well-established analytic and algebraic tools previously used to study lattice paths.
In particular, those are the robustness of D-finite functions
with respect to the Hadamard product,
and the kernel method  \cite{FS09,BM10,BaFl02}.

The N-walks are nondeterministic one-dimensional discrete walks.
We will see that their generating functions require three variables:
one marking the lowest point $\min(w)$ that can be reached by the N-walk $w$,
another one marking the highest point $\max(w)$,
and the last one marking its length $|w|$.
Hence, they are also closely related to two-dimensional lattice paths,
if we interpret $\left(\min(w), \max(w)\right)$ as coordinates in the plane.

		\section{Dyck N-walks}

The step set of classical Dyck paths is $\{-1,1\}$.
The  N-step set of all nonempty subsets is
$$S = \big\{ \{-1\}, \{1\}, \{-1, 1\} \big\},$$
and we call the corresponding N-walks \emph{Dyck N-walks}.
To every step we associate a weight or probability 
$p_{-1}, p_{1}$, and $p_{-1,1}$, respectively.

\begin{example}[Dyck $N$-walks]
Let us consider the Dyck N-walk
$w = (\{1\}, \{-1,1\},\{-1,1\}, \{-1\}).$
The sequence of its reachable points is
$\left(\{0\}, \{1\}, \{0,2\}, \{-1,1,3\}, \{-2,0,2\} \right)$.
There are $4$ classical walks compatible with it:

\begin{center}
\begin{tabular}{|c|c|}
\hline
Classical walk & Geometric realization \\
(sequence of steps) & (ordinates)\\
\hline
$(1,-1,-1,-1)$&$(0,1,0,-1,-2)$\\
$(1,-1,1,-1)$&$(0,1,0,1,0)$\\
$(1,1,-1,-1)$&$(0,1,2,1,0)$\\
$(1,1,1,-1)$&$(0,1,2,3,2)$
\\
\hline
\end{tabular}
\end{center}

\noindent There are two bridges, which happen to be excursions.
Thus, $w$ is an N-bridge and an N-excursion.

\end{example}

The set of reachable points of a Dyck N-walk
or N-meander has the following particular structure.

\begin{lemma} \label{th:dyck_reachable_points}
The reachable points of a Dyck N-walk $w$ are
$
    \left\{\min(w) + 2 i\ |\ 0 \leq \min(w) + 2i \leq \max(w)\right\},
$
where $\min(w)$, $\max(w)$, and the length of $w$ have the same parity.
The same result holds for Dyck N-meanders,
with $\min(w)$ and $\max(w)$ replaced by $\min^+(w)$ and $\max^+(w)$
(see Definition~\ref{def:reachable_points}).
\label{lemma:minmax}
\end{lemma}


We define the generating functions $D(x,y; t)$, $D^+(x,y; t)$,
of Dyck N-walks and Dyck N-meanders as
\begin{align*}
  \sum_{\text{Dyck N-walk $w$}}
   &\bigg( \prod_{s \in w} p_s \bigg)
  x^{\min(w)}
  y^{\max(w)}
  t^{|w|},
  \\
  \sum_{\text{Dyck N-meander $w$}}
  &\bigg( \prod_{s \in w} p_s \bigg)
  x^{\min^+(w)}
  y^{\max^+(w)}
  t^{|w|}.
\end{align*}
Note that by construction these are power series in $t$ with Laurent polynomials in $x$ and $y$, as each of the finitely many N-walks of length $n$ has a finite minimum and maximum reachable point.

\begin{Remark}
One difference to classical lattice paths is the choice of the catalytic variables $x$ and $y$.
Here, they encode the minimum and the maximum reachable points, while in classical problems one chooses to keep track of the coordinates of the endpoint, (see \cite{BaFl02}, for example).
\end{Remark}

\subsection{Dyck N-meanders and N-excursions}

As a direct corollary of Lemma~\ref{lemma:minmax},
all N-bridges and N-excursions have even length.
The total number of Dyck N-bridges and Dyck N-excursions are then, respectively, given by
\[
  [x^{\leq 0} y^{\geq 0} t^{2n}] D(x,y; t)
  \qquad \text{and} \qquad
  D^+(0,1; t),
\]
where the coefficient extraction operator $[t^k]$ is defined as $[t^k] \sum_{n\geq 0} f_n t^n := f_k$
and the nonpositive part extraction operator $[x^{\leq 0}]$ is defined as
$[x^{\leq 0}] \sum_{k \in \Z} g_k x^k := \sum_{k \leq 0} g_k x^k$
(and analogously for $[y^{\geq 0}]$).

\begin{proposition} \label{th:Dyck_Nmeanders}
The generating function of Dyck N-meanders is characterized by the relation
\begin{align*}
  D^+(x,y; t) =
  1 &+
  t \left( p_{-1} x^{-1} y^{-1} + p_1 x y + p_{-1,1} x^{-1} y \right)\\
  & \qquad \times (D^+(x,y; t) - D^+(0,y; t))
  \\ &+
  t \left(p_{-1} x y^{-1} + (p_1 + p_{-1,1}) x y \right)\\
  & \qquad \times (D^+(0,y; t) - D^+(0,0; t))\\
  &+
  t \left( p_1 + p_{-1,1} \right) x y D^+(0,0; t).
\end{align*}
\end{proposition}

\begin{proof}
Applying the symbolic method (see \cite{FS09}),
we translate the following combinatorial characterization
of N-meanders into the claimed equation.
An N-meander is either of length $0$,
or it can be uniquely decomposed into an N-meander $w$
followed by an N-step.
If $\min^+(w)$ is nonzero, then any N-step can be applied.
The generating function of N-meanders
with positive minimum reachable point
is $D^+(x,y; t) - D^+(0,y; t)$.
If $\min^+(w)$ vanishes, but $\max^+(w)$ is nonzero
(those N-meanders have generating function $D^+(0,y; t) - D^+(0,0; t)$),
then an additional N-step $\{-1\}$ increases $\min^+(w)$
(the path ending at~$0$ disappears, and the one ending at~$2$ becomes the minimum) and decreases $\max^+(w)$,
while an additional N-step $\{1\}$ or $\{-1,1\}$ increases both $\min^+(w)$ and $\max^+(w)$.
Finally, if $\min^+(w)$ and $\max^+(w)$ vanish,
which corresponds to the generating function $D^+(0,0; t)$,
then the N-step $\{-1\}$ is forbidden,
and the two other available N-steps both increase $\min^+(w)$ and $\max^+(w)$.
\end{proof}

Let us introduce the \emph{min-max-change polynomial} $S(x,y)$
and the \emph{kernel} $K(x,y)$ as
\begin{align*}
	S(x,y) &:= \frac{p_{-1}}{x y} + p_{1} x y + p_{-1,1} \frac{y}{x},
  \\
  K(x,y)  &:= xy(1- t S(x,y)).
\end{align*}
%
The generating function of Dyck N-walks has now the compact form
$1/(1 - t S(x,y))$.
A key role in the following result on the closed form of Dyck N-meanders
is played by $Y(t)$ and $X(y,t)$, the unique power series solutions satisfying $K(1,Y(t)) = 0$, and $K(X(y,t),y)=0$ which are given by
\begin{align*}
	Y(t) &= \frac{1-\sqrt{1-4p_{-1}(p_{1}+p_{-1,1})t^2}}{2(p_{1}+p_{-1,1})t},
  \\
  X(y,t) &= \frac{1-\sqrt{1-4p_{1}(p_{-1}+p_{-1,1}y^2)t^2}}{2p_{1} y t}.
\end{align*}

\begin{theorem} \label{th:other_Dyck_Nmeanders}
	The generating function $D^+(x,y;t)$ of Dyck $N$-meanders is algebraic of degree $4$, and equal to
    \[
    	 \frac{x-X(y,t)}{1-X(y,t)^2} \frac{y-xY(t)-X(y,t)Y(t)+xyX(y,t)}{xy(1-tS(x,y))}.
    \]
    The generating function of Dyck N-excursions
    is symmetric in $p_{-1}$ and $p_{1}$, and equal to
    \begin{align*}
    	D^+(0,1;t) &= \frac{X(1,t)}{1-X(1,t)^2} \frac{1-X(1,t)Y(t)}{(p_{-1}+p_{-1,1})t}.
    \end{align*}
\end{theorem}

\begin{proofof}{Proof (Sketch)}
Starting from the result of Proposition~\ref{th:other_Dyck_Nmeanders} one first substitutes $x=1$ and finds a closed-form expression for $D^+(0,0;t)$ using the kernel method. After substituting this expression back into the initial equation one applies the kernel method again with respect to $x$ and finds a closed-form solution for $D^+(0,y;t)$. Combining these results one proves the claim.
Finally, using a computer algebra system a short computation using the closed form of Dyck N-excursions shows the symmetry in $p_{-1}$ and $p_{1}$.
\end{proofof}

\begin{Remark}
	It would be desirable to find a combinatorial interpretation of the surprising symmetry in $p_{-1}$ and $p_{1}$ of Dyck N-excursions (which is clear for Dyck N-bridges).
\end{Remark}

 
With this result, we can easily answer the counting problem in which all weights are set equal to one. Thereby we also solve a conjecture in the OEIS\footnote{The on-line encyclopedia of integer sequences: \url{http://oeis.org/A151281}.} on the asymptotic growth.

\begin{corollary}
  For $p_{-1}=p_{1}=p_{-1,1}=1$ the generating function of unweighted Dyck N-meanders is
  \begin{align*}
      D^+(1,1,t)
      &= -\frac{1-4t-\sqrt{1-8t^2}} {4t(1-3t)} 
      \\&= 1 + 2t + 6t^2 + 16t^3 + 48t^4 + \ldots.
  \end{align*}
  The number of unweighted Dyck N-meanders is asymptotically equal to
  \begin{align*}
  	[t^n]D^+(1,1,t) = \frac{3^n}{2} + & \left( 3\sqrt{2}(1+(-1)^n) + 4(1-(-1)^n) \right) \\
    & \times \frac{8^{n/2}}{\sqrt{\pi n^3}} + \bigO\left(\frac{8^{n/2}}{n^{5/2}}\right).
  \end{align*}
  These N-walks are in bijection with walks in the first quadrant $\Z_{\geq 0}^2$ starting at $(0,0)$ and consisting of steps $\{(-1,0),(1,0),(1,1)\}$.
  The counting sequence is given by \OEIS{A151281}. 
  
  For $p_{-1}=p_{1}=p_{-1,1}=1$ the complete generating function of unweighted Dyck N-excursions is
  \begin{align*}
      D^+(0,1,t) &= \frac{1-8t^2-(1-12t^2)\sqrt{1-8t^2}} {8t^2(1-9t^2)} 
      \\&= 1 + 4t^2 + 28t^4 + 224^6 + 1888^8 + \ldots.
  \end{align*}	
  The number $[t^n]D^+(0,1,t)$ of unweighted Dyck N-excursions is asymptotically equal to
  \begin{align*}
  	 (1+(-1)^n)\left( \frac{3^n}{8} + \sqrt{8} \frac{8^{n/2}}{\sqrt{\pi n^3}}  + \bigO\left(\frac{8^{n/2}}{n^{5/2}}\right) \right).
  \end{align*}
\end{corollary}


Finally, we come back to one of the starting questions from the networking motivation.

\begin{theorem}
The probability for a random Dyck N-walk
of length $2n$ to be an N-excursion has for $n \to \infty$ the following asymptotic form where the roles of $p_{-1}$ and $p_{1}$ are interchangeable:
\medskip
\begin{itemize}
\item $\frac{(1-2p_{1})(1-2p_{-1})}{(1-p_{1})(1-p_{-1})} + \bigO\left(\frac{\left(4p_{-1}(1-p_{-1})\right)^n}{n^{3/2}}\right)$
if $0 < p_{1} \leq p_{-1} < \frac{1}{2}$, 
\item $\frac{1-2p_{1}}{(1-p_{1})\sqrt{\pi n}} + \bigO\left(\frac{1}{n^{3/2}}\right)$
if $0 < p_{1} < \frac{1}{2}$ and $p_{-1} = \frac{1}{2}$,
\item $\frac{1}{\sqrt{\pi n^3}} + \bigO\left(\frac{1}{n^{5/2}}\right)$
if $p_{1} = p_{-1} = \frac{1}{2}$,
\item $\bigO\left(\frac{\left(4p_{-1}(1-p_{-1})\right)^n}{n^{3/2}}\right)$
if $0 < p_{1} < \frac{1}{2} < p_{-1} < 1$ and $p_{-1} + p_{1} \leq 1.$
\end{itemize}
\end{theorem}

\begin{proofof}{Proof (Sketch)}
	Starting from the results of Theorem~\ref{th:Dyck_Nmeanders} we perform a singularity analysis~\cite{FS09}. Thereby different regimes need to be considered, leading to the different cases in the result. In the last case the condition guarantees that $p_{-1}$ is closer to $1/2$ than $p_{1}$.
\end{proofof}

\begin{figure*}[ht!]
	\centering
	\subfloat[\scriptsize Short paths (from $0$ to $150$)]{%
    \resizebox{0.47\textwidth}{!}{\input{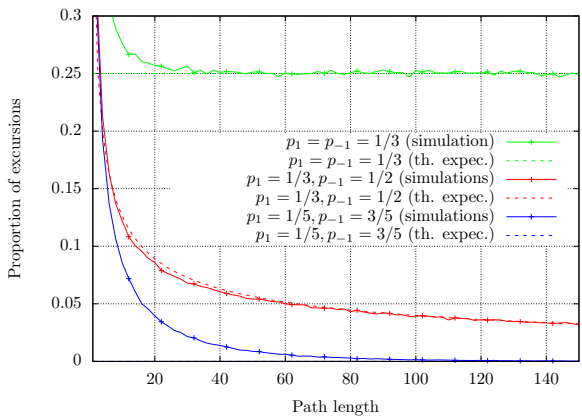}}
		\label{fig:short_dyck}}
        \hfill
	\subfloat[\scriptsize Long paths (from $0$ to $10^4$)]{%
    \resizebox{0.47\textwidth}{!}{\input{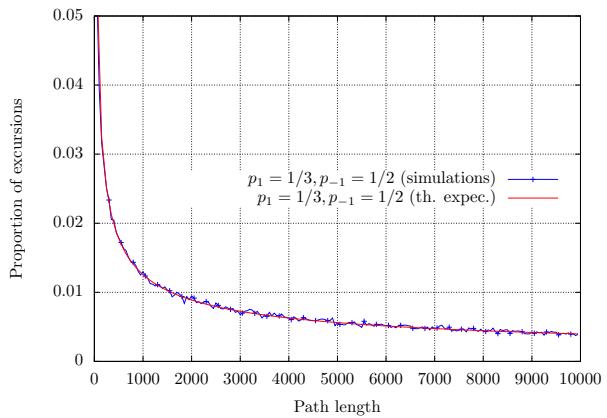}}
		\label{fig:long_dyck}}
    \caption{Comparison of theoretical expectation and averaged simulation (over $10^5$ runs) of the proportion of Dyck N-excursions among Dyck N-walks.}
	\label{fig2} 
\end{figure*}

Note that the (huge) formula for the constant in the last case can be made explicit in terms of $p_{-1}$ and $p_{1}$. However, it is of different shape for $p_{-1} + p_{1} = 1$, and $p_{-1} + p_{1} < 1$. In Figure~\ref{fig2} we compare the theoretical results with simulations for three different probability distributions. These nicely exemplify three of the four possible regimes of convergence.

\subsection{Dyck N-bridges}

We now turn our attention to Dyck N-bridges. Their generating function is defined as
\[
	B(x,y,t) = \sum_{n,k,\ell \geq 0} b_{2n,k,\ell} x^{-k} y^{\ell} t^{2n}.
\]
Recall the following relation with all N-walks (note that bridges have to be of even length):
$
	[t^{2n}]B(x,y,t) = [x^{\leq 0} y^{\geq 0} t^{2n}] D(x,y; t).
$
In the following theorem we will reveal a great contrast to classical walks: nearly all N-walks are N-bridges.

\begin{theorem}
	\label{theo:DyckNBridges}
	The generating function of Dyck N-bridges $B(x,y,t)$ is algebraic of degree $4$. 
    For $\wt{-1}=\wt{1}=\wt{-1,1}=1$ the generating function of unweighted Dyck N-bridges is algebraic of degree $2$:
    \begin{align*}
    	B(1,1,t) &= \frac{1-6t^2}{\sqrt{1-8t^2}(1-9t^2)}\\
      &= 1 + 7t^2 + 63t^4 + 583t^6 + 5407t^8 + \ldots.
    \end{align*}
    The number $[t^n] B(1,1,t)$ of unweighted Dyck N-bridges is asymptotically equal to
    \begin{align*}
    	\frac{1+(-1)^n}{2} \left( 3^n - \frac{2\sqrt{2}}{\sqrt{\pi}} \frac{8^{n/2}}{\sqrt{n}} + \bigO\left(\frac{8^{n/2}}{n^{3/2}}\right) \right).
    \end{align*}
\end{theorem}

\begin{proof}
In order to improve readability we drop the parity condition on $t$ and define
\[
  \Bb(x,y,t) := [x^{\leq 0} y^{\geq 0}] D(x,y; t),
\]
such that
\begin{align} \label{eq:NDyckBridgesInterpretation}
	\begin{aligned}
	\Bb(x,y,t) &= D(x,y; t) - [x^{>0}] D(x,y,t) \\
    		   &  \qquad - [y^{<0}] D(x,y,t).
    \end{aligned}
\end{align}
It is then simple to recover $B(x,y,t)$ from $\Bb(x,y,t)$. 
In words, an N-bridge is an N-walk of even length whose minimum is not strictly positive, nor is its maximum strictly negative\footnote{We thank Mireille Bousquet-M\'elou for suggesting us this approach.}. 
%
%
	
 The change in the $x$- (minimal reachable point) and $y$-coordinate (maximal reachable point) can be conveniently encoded in the min-max-change polynomial
\[
 	S(x,y) = \frac{\wt{-1}}{xy} + \wt{1}xy  + \wt{-1,1}\frac{y}{x}.
\]
Then, the construction can be interpreted as the one of two-dimensional walks of length $n$, starting at $(0,0)$, with the step set $\{(-1,-1),(1,1),(-1,1)\}$, and ending in the fourth quadrant $\{(x,y) : x\geq 0, y\leq 0\}$. A lot is known about these walks, see \eg~\cite{BM10}. By~\eqref{eq:NDyckBridgesInterpretation} it suffices to find the generating functions $F(x,y,t) := [x^{>0}] D(x,y,t)$ 
and $G(x,y,t) := [y^{<0}] D(x,y,t)$ for $2D$-walks ending with a positive abscissa (\resp negative ordinate). 
The theory of 
formal Laurent series with positive coefficients tells us automatically that they are algebraic, 
which implies that the generating function of bridges is algebraic, see \eg~\cite[Section 6]{Gessel80} which also gives further historical references.


Due to the symmetry of the step set we have $F(x,y,t)=G(1/y,1/x,t)$ after additionally interchanging the role of $\wt{-1}$ and $\wt{1}$. 
In order to end the proof it remains to compute the roots of the denominator of $D(x,y,t)$ and perform a partial fraction decomposition. 
%
\end{proof}


After this detailed discussion of nondeterministic walks derived from Dyck paths, we turn to the probably next most classical lattice paths: Motzkin paths.

		\section{Motzkin N-walks}

The step set of classical Motzkin paths is $\{-1,0,1\}$.
The  N-step set of all nonempty subsets is
\[
	S = \big\{\{-1\}, \{0\}, \{1\}, \{-1,0\}, \{-1,1\}, \{0,1\}, \{-1,0,1\} \big\},
\]
and we call the corresponding N-walks \emph{Motzkin N-walks}.
%
A Motzkin N-walk $w$ is said to be 
\begin{itemize}
  \item of type $1$ if $\reach(w)$ is equal to $\left\{\min(w), \min(w) + 2, \min(w) + 4, \ldots, \max(w)\right\}$,
  \item of type $2$ if $\reach(w)$ is equal to $\left\{\min(w), \min(w) + 1, \min(w) + 2, \ldots, \max(w)\right\}$
    and $\max(w) - \min(w) \geq 1$.
\end{itemize}
The following proposition explains how these two types are sufficient to characterize the structure of Motzkin N-walks.

\begin{proposition} \label{th:motzkin_structure}
A Motzkin N-walk $w$ is of type $1$ if and only if
it is constructed only from the N-steps  $\{-1\}$, $\{0\}$, $\{1\}$, and $\{-1,1\}$.
Otherwise, it is of type $2$.
\end{proposition}

\begin{proofof}{Proof (Sketch)}
The proof is based on a recurrence and a case-by-case analysis on the number and type of N-steps.
\end{proofof}

The set of Motzkin N-walks of type $1$ (\resp $2$) is denoted by $M_1$ (\resp $M_2$),
and their generating functions are defined as
\begin{align*}
	M_1(x,y; t) &= \sum_{w \in M_1} x^{\min(w)} y^{\max(w)} t^{|w|},\\
	M_2(x,y; t) &= \sum_{w \in M_2} x^{\min(w)} y^{\max(w) - 1} t^{|w|}.
\end{align*}

\begin{theorem}
The generating functions of Motzkin N-walks of type $1$ and $2$ are rational.
The generating function of Motzkin N-bridges is algebraic.
\end{theorem}

\begin{proof}
The first statement is a direct corollary of the previous proposition due to a simple sequence construction.
An N-bridge $w$ of type $1$ is an $M_1$ N-walk that satisfies $\min(w) \leq 0$, $\max(w) \geq 0$, and $\min(w)$ is even.
Note that in this case this property is not equivalent to an even number of steps.
An N-bridge $w$ of type $2$ is an $M_2$ N-walk that satisfies $\min(w) \leq 0$ and $\max(w) \geq 0$.
Thus, the generating function of Motzkin N-bridges is equal to
\[
    [x^{\leq 0} y^{\geq 0}] \left(\frac{M_1(x,y; t) + M_1(-x,y; t)}{2} + M_2(x,y; t)
    \right).
\]
Since the generating functions of $M_1$ and $M_2$ are rational,
according to \cite[Proposition~1]{BM10} (see also \cite{Lipshitz88}), the generating function of N-bridges is D-finite. 
Yet the generating function is even algebraic, which can be proved similarly as done the proof of Theorem~\ref{theo:DyckNBridges}.
\end{proof}

\begin{Remark}
Using a computer algebra system it is easy to get closed-form solutions and asymptotics for specific values of the weights. 
We do not give these closed forms, as they are quite large and do not shed new light on the problem.
It is however interesting to compute the asymptotic proportion of N-bridges among all N-walks. For example, when all weights are set to $1$, it is equal to
\[
	1 - \sqrt{\frac{3}{\pi}} \frac{(6/7)^n}{\sqrt{n}} + \bigO\left(\frac{(6/7)^n}{n^{3/2}}\right).
\]
Hence, nearly all N-walks are N-bridges.
\end{Remark}

We now turn to the analysis of Motzkin N-meanders and N-excursions.

\begin{theorem} \label{th:motzkin_N_meanders}
The generating functions of Motzkin N-meanders and N-excursions are algebraic.
\end{theorem}

\begin{proof}
Without loss of generality we perform all computations here with all weights $p_i=1$. 
Let $M_1^+$ and $M_2^+$ denote the Motzkin N-meanders of type $1$ and $2$.
Their generating functions are
\begin{align*}
	M_1^+(x,y; t) &= \sum_{w \in M_1^+} x^{\min^+(w)} y^{\max^+(w)} t^{|w|},\\
	M_2^+(x,y; t) &= \sum_{w \in M_2^+} x^{\min^+(w)} y^{\max^+(w) - 1} t^{|w|}.
\end{align*}
Let also $M^+(x,y; t)$ denote the column vector $(M_1^+(x,y; t), M_2^+(x,y; t))$.
An N-meander is either empty -- in which case, it is of type $1$ --
or it is an N-meander $w$ followed by an N-step $s$.
The type of $w \cdot s$ depends on the type of $w$,
the N-step $s$, as well as on the case if $\min^+(w) = 0$ or if $\max^+(w) = 0$.
Specifically,
\begin{itemize}
\item
when $w$ has type $1$,
then $w \cdot s$ has type $1$ if $s \in \{\{-1\}, \{0\}, \{1\}, \{-1,1\}\}$, otherwise it has type $2$,
\item
when $w$ has type $2$ and $\max^+(w) > 1$ then $w \cdot s$ has type $2$ for any $s$,
\item
when $w$ has type $2$ and $\max^+(w) = 1$ (i.e.~the reachable points are $\{0,1\}$) then $w \cdot s$ has type $1$ if $s = \{-1\}$, and type $2$ otherwise.
\end{itemize}
Applying the Symbolic Method \cite{FS09} and the same reasoning as in the proof of Proposition~\ref{th:Dyck_Nmeanders},
we obtain the following system of equations
characterizing the generating functions from the vector $M^+(x,y; t)$
\begin{align*}
	M^+(x,y; t) &= 
    e_1 + t \Big( 
    A(x,y) (M^+(x,y; t) - M^+(0,y; t)) \\
    &\quad+ B(x,y) (M^+(0,y; t) - M^+(0,0; t)) 
    \\ & \quad+ C(x,y) M^+(0,0; t) \Big),
\end{align*}
where $e_1$ is the column vector $(1,0)$,
and $A(x,y)$, $B(x,y)$, $C(x,y)$ are two-by-two matrices with Laurent polynomials in $x$ and $y$ given in Figure~\ref{fig:matrices}.
Observe that the first two matrices are upper-triangular.
\begin{figure*}
\begin{align*}
	A(x,y) &= {\small \begin{pmatrix} x^{-1} y^{-1} + 1 + x y + x^{-1} y & 0\\ x^{-1} y^{-1} + 1 + x^{-1} & x^{-1} y^{-1} + 1 + x y + x^{-1} + y + 2 x^{-1} y \end{pmatrix}},\\
    B(x,y) &= {\small \begin{pmatrix} x y^{-1} + 1 + 2 x y & 0\\ y^{-1} + 2 & y^{-1} + 2 + x y + 3 y \end{pmatrix}}, \qquad
    C(x,y) = {\small \begin{pmatrix} 2 + 2 x y & 1\\ 2 & x y + 2 + 3 y \end{pmatrix}}.
\end{align*}
    \caption{Matrices involved in the proof of Theorem~\ref{th:motzkin_N_meanders}.}
		\label{fig:matrices}
\end{figure*}
This equation is rearranged into
\begin{align} 
	&\left(\id - t A(x,y)\right) M^+(x,y; t) = \nonumber
  \\&
    e_1 - t \left(A(x,y) - B(x,y)\right) M^+(0,y; t) \label{eq:motzkin_n_meanders}
  \\& \quad - t \left(B(x,y) - C(x,y)\right) M^+(0,0; t). \nonumber
\end{align}
Next, we apply the kernel method (see \eg \cite{BaFl02} and \cite{AsBaBaGi18})
successively on $x$ and $y$ in a two phases
to compute the generating function $M^+(x,y; t)$ of Motzkin N-meanders.
The small roots in the variable $x$ of the equations
\begin{align*}
	1 - t A_{0,0}(x,y) &= 0, 
    \\
    1 - t A_{1,1}(x,y) &= 0,
\end{align*}
are denoted by $X_1(y,t)$ and $X_2(y,t)$, and are equal to
\begin{gather*}
	\frac{1 - t - \sqrt{1 - 4 t^2 y^2 - 3 t^2 - 2 t}}{2 t y},
  \\
  \frac{1 - t (y + 1) - \sqrt{1 - 7 t^2 y^2 - 2 t^2 y - 3 t^2 - 2 t y - 2 t}}{2 t y}.
\end{gather*}
We then define the row vectors
\begin{align*}
	u_1 &= (1, 0),
  \\
  u_2(y,t) &= \left(t A_{1,0}(X_2(y,t), y), 1 - t A_{0,0}(X_2(y,t), y)\right),
\end{align*}
so that the left-hand side of Equation~\eqref{eq:motzkin_n_meanders}
vanishes both when evaluated at $x = X_1(y,t)$ and left-multiplied by $u_1$,
and also when evaluated at $x = X_2(y,t)$ and left-multiplied by $u_2(y,t)$.
Combining the corresponding two right-hand sides, we obtain a new two-by-two system of linear equations
\begin{equation} \label{eq:motzkin_n_meanders_second_kernel}
	t D(y,t) M^+(0,y; t) = f(y,t) -  E(y,t) M^+(0,0; t),
\end{equation}
where the vector $f(y,t)$ of size $2$ has its first element equal to $1$,
and its second element equal to
\[
	\frac{{\left(t y + t - \sqrt{-7 \, t^{2} y^{2} - 3 \, t^{2} - 2 \, {\left(t^{2} + t\right)} y - 2 \, t + 1} + 1\right)} t}
  {1 - t y - t - \sqrt{-7 \, t^{2} y^{2} - 3 \, t^{2} - 2 \, {\left(t^{2} + t\right)} y - 2 \, t + 1}},
\]
and the two-by-two matrices $D(y,t)$ and $E(y,t)$ are two large to be shown here.
Again, the matrix $D(y,t)$ is upper-triangular.
We now define
\begin{align*}
	Y_1(t) &= \frac{t-1+\sqrt{-7 t^2-2 t+1}}{4 t},
  \\
  Y_2(t) &= \frac{1 - 2 t - \sqrt{-12 t^2-4 t+1}}{8 t},
\end{align*}
and the row vectors
\begin{align*}
	v_1 &= (1,0),
    \\
    v_2(t) &= (- D_{1,0}(Y_2(t),t), D_{0,0}(Y_2(t), t)),
\end{align*}
to ensure that $Y_1(t)$ and $Y_2(t)$ have series expansions at the origin,
and that the left-hand side of Equation~\eqref{eq:motzkin_n_meanders_second_kernel}
vanishes both when evaluated at $y = Y_1(t)$ and left-multiplied by $v_1$,
and also when evaluated at $y = Y_2(t)$ and left-multiplied by $v_2(t)$.
The corresponding two right-hand side are combined to form a new two-by-two system of equations
\[
	h(t) = t F(t) M^+(0,0; t),
\]
where the column vector $h(t)$ and the matrix $F(t)$ are too large to be shown here.
The matrix $F(t)$ is invertible, so the generating function of Motzkin N-meanders with maximum reachable point $0$ is equal to
\[
	M^+(0,0; t) = \frac{1}{t} F(t)^{-1} h(t).
\]
This expression is injected in Equation~\eqref{eq:motzkin_n_meanders_second_kernel}
to express the generating function of Motzkin N-meanders with minimum reachable point $0$
\[
	M^+(0,y; t) = \frac{1}{t} D(y,t)^{-1} (f(y,t) -  E(y,t) M^+(0,0; t)).
\]
Finally, this expression is injected in Equation~\eqref{eq:motzkin_n_meanders}
to express the generating function of Motzkin N-meanders
\begin{align*}
	M^+(x,y; t) =&
	\left(\id - t A(x,y)\right)^{-1} 
  \\ & \times (e_1 - t \left(A(x,y) - B(x,y)\right) M^+(0,y; t) 
  \\ & \qquad - t \left(B(x,y) - C(x,y)\right) M^+(0,0; t)).
\end{align*}
The generating function of N-meanders and N-excursions are then, respectively,
$M_1^+(1,1; t) + M_2^+(1,1; t)$ and $M_1^+(0,1; t) + M_2^+(0,1; t)$.
\end{proof}

\begin{Remark}
As before we can use a computer algebra system to get numeric results. After tedious computations one gets that for all $p_i$'s equal to $1$ the generating function of N-meanders is algebraic of degree $2$ and given by
\begin{align*}
	&\frac{10t-1+\sqrt{(1+2t)(1-6t)}}{8t(1-7t)}.
\end{align*}
The total number of N-meanders is asymptotically equal to
\begin{align*}
	\frac{3}{4}7^n + \frac{3\sqrt{3}}{2 \sqrt{\pi}} \frac{6^n}{\sqrt{n^3}} + \bigO\left(\frac{6^n}{n^{5/2}}\right).
\end{align*}
The generating function of N-excursions is algebraic of degree $4$.
Their asymptotic number is
\begin{align*}
	\frac{9}{16}7^n - \gamma \frac{6^n}{\sqrt{\pi n^3}} + \bigO\left(\frac{6^n}{n^{5/2}}\right), \quad 
\end{align*}
where $\gamma \approx 0.6183$ is the positive real solution of $1024\gamma^4-8019\gamma^2+2916=0$. 
This means that for large $n$ approximately $75\%$ of all N-walks are N-meanders and $56.25 \%$ of all N-walks are N-excursions.
\end{Remark}

\begin{figure}
\begin{center}
\includegraphics[scale=0.35]{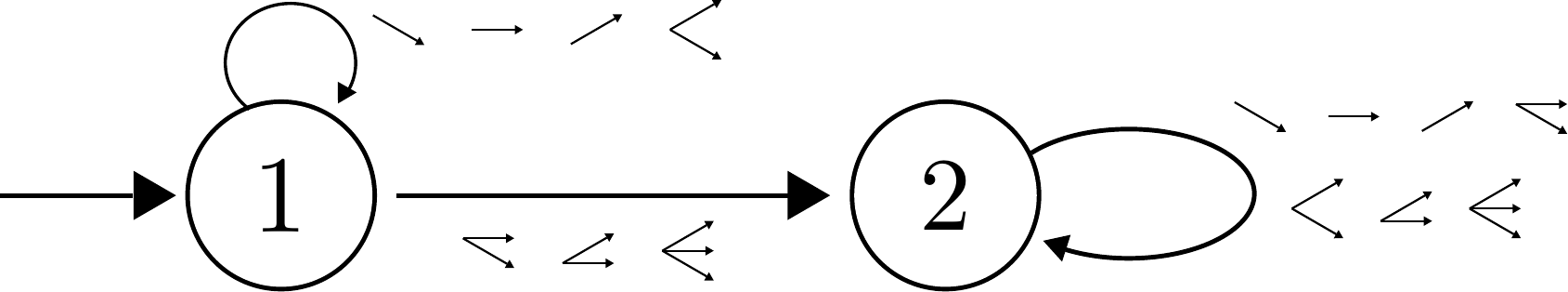}
\caption{The automaton representing the structure of reachable points of Motzkin N-walks.
The types from Theorem~\ref{th:general_structure} corresponding to vertex $1$ are
$A_1 = \{0\}$,
$B_1 = \{1\}$,
$C_1 = \emptyset$,
and for vertex $2$, we have
$A_2 = \{0\}$,
$B_2 = \{0\}$,
$C_2 = \{0\}$.}
\label{fig:motzkin_structure}
\end{center}
\end{figure}

		\section{N-bridges with general N-steps}

The main result of this section is

\begin{theorem} \label{th:general_bridges}
For any N-step set $S$, the generating function of N-bridges is algebraic.
\end{theorem}

A method for computing this generating function is provided by the proof,
in Subsection~\ref{sec:proof_general_bridges}.
In order to establish this result,
we first derive Proposition~\ref{th:general_structure},
which provides a description of the set of reachable points
for N-walks on a given N-step set.
It is proven in Section~\ref{sec:proof_general_structure}.

Given nonnegative integer sets $A$, $B$, $C$,
an N-walk $w$ is \emph{of type $(A, B, C)$}
when an integer $r$ is reachable if and only if
$\max(w) - \min(w) \geq \max(A) + \max(C)$,
and at least one of the following conditions holds
\begin{itemize}
\item
$r - \min(w)$ belongs to $A$,
\item
$\max(w) - r$ belongs to $C$,
\item
$r - \min(w) \geq \max(A)$,
$\max(w) - r \geq \max(C)$,
and $r - \min(w) - \max(A) - 1 \mod (\max(B) + 1)$ belongs to $B$,
\end{itemize}
with the convention $\max(\emptyset) = 0$.
The set of N-walks of type $(A,B,C)$ is denoted by $W_{A,B,C}$.

\begin{proposition} \label{th:general_structure}
Given an N-step set $S$, there is a finite set of types $(A_j,B_j,C_j)_{1 \leq j \leq m}$
such that the set of all N-walks on $S$ is equal to the disjoint union
$
	\biguplus_{j=1}^m W_{A_j,B_j,C_j}.
$
Furthermore, if we consider the N-walks as words on the alphabet $S$,
there are nonempty subsets $(T_{i,j})_{1 \leq i \leq j \leq m}$ of $S$
such that for all $1 \leq j \leq m$,
the grammar characterizing the N-walks of type $(A_j,B_j,C_j)$ is
\begin{align*}
	W_{A_j,B_j,C_j} =& \text{ $($empty N-walk if $j=1)$} \\& \quad + \sum_{i=1}^j W_{A_i,B_i,C_i} \sum_{s \in T_{i,j}} s.
\end{align*}
\end{proposition}

Figure~\ref{fig:motzkin_structure} illustrates the previous proposition
on the example of Motzkin N-walks.

	\subsection{Proof of Proposition~\ref{th:general_structure}} \label{sec:proof_general_structure}

Consider an N-walk $w$ and the N-walk $w \cdot s$ obtained by adding the N-step $s$ to $w$.
In this section, we will use the observation that the set of reachable points of $w$ and $w \cdot s$
are linked by the relation
\[
	\reach(w \cdot s) = \bigcup_{h \in s} \left\{r + h\ |\ r \in \reach(w)\right\}.
\]
Hence, the set of reachable points of an N-walk does not depend on the order of its N-steps.
We start the proof with a description of reachable points as solutions of a linear equation.

\begin{lemma} \label{th:reachable_linear_equation}
There exists an N-walk on the N-step set $S$ that reaches the point $r$
if and only if the following equation has a solution
\[
	\sum_{s \in S} \sum_{h \in s} h\, x_{s,h} = r,
    \quad
    \forall (s,h),\ x_{s,h} \in \naturals.
\]
Furthermore, any N-walk that contains exactly $\sum_{h \in s} x_{s,h}$
occurrences of the N-step $s$ reaches $r$.
\end{lemma}

\begin{proof}
By definition, if the N-walk $w = s_1 \cdot s_2 \cdots s_{|w|}$ on the N-step set $S$ reaches the point $r$,
then there exists a sequence of integers $(h_i)_{1 \leq i \leq |w|}$
such that for all $i$, we have $h_i \in s_i$,
and
\[
	\sum_{i=1}^{|w|} h_i = r.
\]
Let $x_{h,s}$ denote the number if values of $1 \leq i \leq |w|$
such that $(h, s) = (h_i, s_i)$, then the previous equation becomes
\[
	\sum_{s \in S} \sum_{h \in s} h\, x_{s,h} = r.
\]
\end{proof}

The previous lemma translates the study of reachable points
into the realm of numerical semigroups.
Using the tools of this field (Schur's Theorem and Fr\"obenius number \cite{Al09}),
we obtain the following lemma,
that specializes Proposition~\ref{th:general_structure}
to N-walks containing sufficiently many occurrences of each N-step.

\begin{lemma} \label{th:particular_case_general_structure}
Let $p_S$ denote the gcd of the N-steps from $S$, shifted so that their minimum is at $0$
\[
	p_S = \gcd \bigg( \bigcup_{s \in S} \{h - \min(s)\ |\ h \in s\} \bigg).
\]
For any N-step set $S$, there exist an integer $m_S$ and two nonnegative integer sets $A$ and $C$
such that any N-walk $w$ on $S$ that contains at least $m_S$ occurrences
of each N-step is in $W_{A,\{p_S\},C}$.
\end{lemma}

\begin{proof}
Given an N-step set $S$, the normalized version
of the N-step $s$ is defined as $\{(h - \min(s))/p_S\ |\ h \in s\}$.
The normalized version of $S$ is then the set of its normalized N-steps.
If Lemma~\ref{th:particular_case_general_structure} holds for normalized N-step sets,
it also holds in the general case.
Thus, without loss of generality, we assume $S$ to be normalized.
In particular, all its N-steps have minimum $0$,
so the smallest reachable point is always $0$.
According to Schur's Theorem, there is an integer $f$, called the \emph{Fr\"obenius number}\footnote{Computing the Fr\"obenius number is $\mathsf{NP}$-hard under Turing reduction if the number of integers $n = |\cup_{s \in S} s|$ is arbitrary~\cite{ramirez1996complexity}. It is an open problem whether it is also $\mathsf{NP}$-hard under Karp reduction. If $n$ is fixed, there is a polynomial algorithm~\cite{kannan1992lattice} to compute the Fr\"obenius number but it is unpractical as its complexity is in $O\left((\log m)^{n^{O(n)}}\right)$ where $m=\max_{s \in S} \max(s)$. However, there are algorithms that perform very well in practice\cite{beihoffer2005faster,einstein2007frobenius}.},
such that for any $r > f$, the equation from Lemma~\ref{th:reachable_linear_equation} has a solution.
Let $w_r$ denote an N-walk reaching $r$, and $|w|_s$ the number of occurrences of the N-step $s$ in $w$.
Let us define the integers $\ell$, $m_S^{(0)}$ and $b$ as
\begin{align*}
	\ell &= \max_{s \in S} \max(s),\\
	m_S^{(0)} &= \max_{s \in S, f < r \leq f + \ell} |w_r|_s,\\
	c &= \max_{\forall s \in S,\ |w|_s = m_S^{(0)}} (\max(w) - f - \ell).
\end{align*}
Those three integers have the following meanings:
\begin{itemize}
\item
$\ell$ is the maximum height of any N-step from $S$,
\item
any N-walk containing at least $m_S^{(0)}$ occurrences of each N-step
reaches all integers from $[f + 1, f + \ell]$,
\item
let $W^{= m_S^{(0)}}$ denote the set of N-walks
that contain exactly $m_S^{(0)}$ occurrences of each N-step,
then for any such N-walk, the distance between the maximal reachable point and $f + \ell$ is at most $c$.
\end{itemize}
Since any N-step has minimum $0$, and maximum at most $\ell$,
adding an N-step $s$ to an N-walk $w$ from $W^{= m_S^{(0)}}$
produces an N-walk $w \cdot s$ which reaches all the points from $f+1$ to $\max(w \cdot s) - c$.
By recurrence, for any N-walk $w$ that contains at least $m_S^{(0)}$ occurrences of each N-step,
all points from $f+1$ to $\max(w) - c$ are reachable.
Since all N-steps contain $0$, we have
\[
	\reach(w) \subset \reach(w \cdot s).
\]
Let $w_{= m}$ denote an N-walk that contains exactly $m$ occurrences of each N-step.
Then $([0, f] \cap \reach(w_{= m}))_m$ is an increasing (for the inclusion) sequence of sets included in $[0,f]$.
Thus, it reaches for some finite integer $m = n$ its limit $A$.
We set $m_S^{(1)} = \max(m_S^{(0)}, n)$.
Any N-walk $w$ containing at least $m_S^{(1)}$ occurrences of each N-step satisfies
\begin{align*}
	& \reach(w) \cap [0,f] = A, \\ & [f + 1, \max(w) - c] \subset \reach(w).
\end{align*}
Finally, let us define the symmetric of an N-step $s$ as the N-step $\{\max(s) - h\ |\ h \in s\}$,
and the symmetric of an N-step set as the set of its symmetric N-steps.
Applying the previous proof to the symmetric of $S$, we obtain the existence
of an integer $m_S^{(2)}$, and integer $f'$ and a set $C$ such that
for any N-walk $w$ containing at least $m_S^{(2)}$ occurrences of each N-step,
we have
\[
	\reach(w) \cap [\max(w) - f', \max(w)] = C.
\]
Defining the integer $m_S$ as $\max(m_S^{(1)}, m_S^{(2)})$ finishes the proof.
\end{proof}

We can finally provide the proof of Proposition~\ref{th:general_structure}.

\begin{proofof}{Proof of Proposition~\ref{th:general_structure}}
Given an N-step set $S$, we set	$q_S = \max_{T \subset S} m_T$,
where $m_T$ has been defined in Lemma~\ref{th:particular_case_general_structure}.
Let $|w|_s$ denote the number of occurrences of the N-step $s$ in the N-walk $w$.
Consider the partial order on N-walks such that $w \leq w'$
if and only if any N-step that has less than $q_S$ occurrences in $w'$
has at least as many occurrences in $w'$ as in $w$.
Let us also define the equivalence relation $w \sim w'$
when the N-steps that occurs at least $q_S$ times in $w$ and $w'$ are the same,
and the other N-steps have the same number of occurrences in both N-walks
\begin{align*}
	w \leq w'  \Leftrightarrow & \left( \forall s \in S, |w'|_s < q_S \Rightarrow |w|_s \leq |w'|_s \right),\\
    w \sim w' \Leftrightarrow & \big( \forall s \in S, (|w|_s \geq q_s \Rightarrow |w'|_s \geq q_s) \\
    & \quad \text{~\, and } (|w|_s < q_s \Rightarrow |w|_s = |w'|_s) \big).
\end{align*}
When the set of all N-walks is factored by the ``$\sim$'' relation,
we obtain a finite number of disjoint subsets,
on which the ``$\leq$'' partial order induces a lattice structure.
In the next paragraph, we will prove that each of those subsets $V$ corresponds to a type,
such that there are finite nonnegative integer sets $A$ $B$, $C$ such that $V = W_{A,B,C}$.
This will conclude the proof of the proposition, as the lattice structure
ensures the grammar characterization stated in the second part of the proposition.

First, observe that two N-walks
that contain the same N-steps with the same multiplicities
reach the same set of points.
Consider an element $V$ of the lattice.
By definition, if an N-walk from $V$ contains
less than $q_S$ occurrences of each N-step,
then all N-walks in $V$ contain, for each N-step, the same number of occurrences,
and thus have the same set $R$ of reachable points.
We then define $A$ as the set $R$ shifted by $\min(R)$,
and obtain $V = W_{A,\emptyset,\emptyset}$.
Otherwise, let $T \subset S$ denote the set of N-steps
that occur in the N-walks from $V$ at least $q_s$ times.
Let also $v$ denote any N-walk
with exactly $|w|_s$ occurrences of each N-step $s$ from $S \setminus T$,
and no other N-step.
Let $W_{\geq m_T}$ denote the set of N-walks
that contain at least $m_T$ occurrences of each N-step from $T$,
and no other N-step.
Since $q_S \geq m_T$, for any N-walk $w$ from $V$,
there is an N-walk $w'$ from $W_{\geq m_T}$ such that
the reachable points of $w$ are the same as for $w' \cdot v$,
the concatenation of $w'$ and $v$.
Since $w'$ belongs to $W_{\geq m_T}$, according to Lemma~\ref{th:particular_case_general_structure},
there are integer sets $A'$, $C'$ and an integer $p_T$ such that $w'$ is in $W_{A',\{p_T\},C'}$.
Adding the N-steps from $v$ to $w'$ changes the set of reachable points,
and we obtain integer sets $A$, $B$, $C$ such that $w' \cdot v$ belongs to $W_{A,B,C}$.
\end{proofof}

	\subsection{Proof of Theorem~\ref{th:general_bridges}} \label{sec:proof_general_bridges}

The main idea of the proof is the repeated use of closure properties of algebraic functions, see~\cite{FS09}. 
Let $W_{A,B,C}(x,y; t) = \sum_{w \in W_{A,B,C}} x^{\min(w)} y^{\max(w)} t^{|w|}$
and $\bridge_{A,B,C}(t)$
denote the generating functions of N-walks and N-bridges of type $(A,B,C)$.
The first part of Proposition~\ref{th:general_structure} implies
\[
	\bridge(t) = \sum_{j=1}^m \bridge_{A_j,B_j,C_j}(t).
\]
Hence, the proof 
is complete
once established that each $\bridge_{A_j,B_j,C_j}(t)$ is algebraic.
The grammar characterization from the second part of Proposition~\ref{th:general_structure}
is translated into the following system of equations:
for $j$ from $1$ to $m$, the generating function $W_{A_j,B_j,C_j}(x,y; t)$ is equal to
\[
    \one_{j=1} +
    t \sum_{i=1}^j W_{A_i,B_i,C_i}(x,y; t) 
    \sum_{s \in T_{i,j}} x^{\min(s)} y^{\max(s)}.
\]
Solving this system, we obtain a rational expression for each $W_{A_j,B_j,C_j}(x,y; t)$,
because the sets $T_{i,j}$ are nonempty.
In the following, we consider some $1 \leq j \leq m$, set $(A,B,C) = (A_j,B_j,C_j)$,
and prove that $\bridge_{A,B,C}(t)$ is algebraic.
We assume that those three sets are nonempty,
the other cases being similar.

An N-walk where the minimal reachable point is positive
or the maximal reachable point is negative cannot be a bridge.
Thus, we distinguish three kinds of N-walks that have the potential to be bridges:
\begin{itemize}
\item
$w \in W_{A,B,C}^{\I}$ when $-\max(A) \leq \min(w) \leq 0$,
\item
$w \in W_{A,B,C}^{\II}$ when $\min(w) < -\max(A)$ and $\max(C) < \max(w)$,
\item
$w \in W_{A,B,C}^{\III}$ when $0 \leq \max(w) \leq \max(C)$.
\end{itemize}
The corresponding generating functions are expressed
as sums of positive parts in $x$ and $y$ of rational function in $t$
with Laurent polynomials in $x$ and $y$ coefficients
\begin{align*}
	W_{A,B,C}^{\I}(x,y; t) &= 
    	W_{A,B,C}(x,y; t) 
  \\ & \qquad - [x^{< -\max(A)}] W_{A,B,C}(x,y; t) 
  \\ & \qquad - [x^{> 0}] W_{A,B,C}(x,y; t),\\
    W_{A,B,C}^{\II}(x,y; t) &= 
    	W_{A,B,C}(x,y; t) 
  \\ & \qquad - [x^{\geq - \max(A)}] W_{A,B,C}(x,y; t) 
  \\ & \qquad - [y^{\leq \max(C)}] W_{A,B,C}(x,y; t),\\
    W_{A,B,C}^{\III}(x,y; t) &= 
    	W_{A,B,C}(x,y; t) 
  \\ & \qquad - [y^{< 0}] W_{A,B,C}(x,y; t) 
  \\ & \qquad - [y^{> \max(C)}] W_{A,B,C}(x,y; t).
\end{align*}
The set of bridges from $W_{A,B,C}^{\I}$ is denoted by $\bridge_{A,B,C}^{\I}$,
and the same holds for $\II$ and $\III$.
By definition of the type, we have
\begin{itemize}
\item
$w \in \bridge_{A,B,C}^{\I}$ if and only if $-\min(w) \in A$,
\item
$w \in \bridge_{A,B,C}^{\II}$ if and only if $-\min(w) - \max(A) - 1 \mod \max(B) + 1 \in B$,
\item
$w \in \bridge_{A,B,C}^{\III}$ if and only if $\max(w) \in C$.
\end{itemize}
Those characterizations imply
\begin{align*}
	\bridge_{A,B,C}^{\I}(t) &= \sum_{a \in A} [x^{a}] W_{A,B,C}^{\I}(x^{-1},1; t),
    \\
	\bridge_{A,B,C}^{\III}(t) &= \sum_{c \in C} [y^{c}] W_{A,B,C}^{\III}(1,y; t),
\end{align*}
which are algebraic functions because the set $A$ is finite, and $W_{A,B,C}^{\I}(x^{-1},1; t)$
is an algebraic function analytic in $x$ and $t$ at the origin, so
\[
	[x^{a}] W_{A,B,C}^{\I}(x^{-1},1; t) = \frac{d^a}{d x^a} W_{A,B,C}^{\I}(x^{-1},1; t)_{|x = 0}
\]
(the same reasoning applies to $\III$).
Using the classical relation,
$
	\frac{1}{p} \sum_{k=0}^{p-1} F(e^{2 i \pi k / p}) = \sum_{p|n} [z^n] F(z),
$
valid for any series $F(z)$ and $p > 0$
we obtain that $\bridge_{A,B,C}^{\II}(t)$, equal to
\[
	 \sum_{b \in B} \sum_{(\max(B) + 1) | (n - \max(A) - 1 - b)} [x^n] W_{A,B,C}^{\II}(x^{-1},1; t),
\]
is algebraic as well.
We conclude that the generating function of all N-bridges is algebraic, since
\begin{align*}
	\bridge(t) &= \sum_{j=1}^m \bridge_{A_j,B_j,C_j}(t)
  \\ &=
    \sum_{j=1}^m \bridge_{A_j,B_j,C_j}^{\I}(t) + \bridge_{A_j,B_j,C_j}^{\II}(t) 
  \\ & \qquad + \bridge_{A_j,B_j,C_j}^{\III}(t).
\end{align*}

\section{Conclusion}

In this paper we introduced nondeterministic lattice paths and solved the asymptotic counting problem for such walks of the Dyck and Motzkin type.
The strength of our approach relies on the methods of analytic combinatorics, which allowed us to derive not only the asymptotic main terms but also lower order terms (to any order if needed). 
Furthermore, we showed that for a general step set the generating function of bridges is algebraic. 
In the long version of this work we will extend this setting to excursions and meanders with general N-steps.

The method of choice is the well-established kernel method.
We extended it to a two-phase approach in order to deal with two catalytic variables. 

Additionally to the mathematically interesting model, our nondeterminstic lattice paths have applications in the encapsulation and decapsulation of protocols over networks. 
In the long version of this work we want to further explore this interesting bridge between combinatorics and networking.

\medskip

{\bf Acknowledgements:}
\label{sec:ack}
The authors would like to thank Mireille Bousquet-M\'elou for her useful comments and suggestions. We also thank the three referees for their feedback.
The second author was partially supported by the  H\'ERA project, funded by The French National Research Agency. Grant no.: ANR-18-CE25-0002. The third author was supported by the Exzellenzstipendium of the Austrian Federal Ministry of Education, Science and Research and the Erwin Schr{\"o}dinger Fellowship
of the Austrian Science Fund (FWF):~J~4162-N35.

\end{document}